\documentclass[letterpaper, 10 pt, conference]{ieeeconf}  

\IEEEoverridecommandlockouts                              

\title{\LARGE \bf
Distributed Optimization via Energy Conservation Laws in Dilated Coordinates
}

\usepackage[american]{babel}
\usepackage{mathtools} 
\usepackage{booktabs} 
\usepackage{graphicx}
\usepackage{subcaption}
\usepackage{tcolorbox}
\usepackage{defs}

\author{Kushal Chakrabarti${^1}$ and Mayank Baranwal$^{1,2}$
\thanks{$^{1}$KC and MB are with the Data \& Decision Sciences division of Tata Consultancy Services Research, Mumbai 400607, India.}
\thanks{$^{2}$MB is also with the Faculty of Systems \& Control group, Indian Institute of Technology Bombay, Mumbai 400076, India.}%
}

\begin{document}

\maketitle
\thispagestyle{empty}
\pagestyle{empty}

\begin{abstract}
Continuous-time models can reveal accelerated structures in distributed optimization, but their rates need not survive direct discretization. We introduce a second-order primal--dual flow for smooth convex distributed optimization and construct an exactly conserved energy that yields an $\mathcal O(t^{-2})$ rate for both the aggregate objective gap and the squared consensus error. We then prove a horizon-wise $\Omega(k^{-1})$ lower bound for a broad class of single-loop finite-memory primal--dual discretizations, ruling out a $\mathcal O(k^{-2})$ aggregate-objective guarantee within this class. Motivated by this barrier, we develop a double-loop method that combines finite-step polynomial consensus with an accelerated outer update. It uses one gradient evaluation and at most $m-1$ communication rounds per outer iteration, $m$ being the number of agents, maintains exact consensus and achieves an $\mathcal O(k^{-2})$ aggregate-objective rate. Numerical comparisons with representative distributed methods support the theory and quantify the communication cost of acceleration.
\end{abstract}

\section{Introduction}
\label{sec:Intro}

Advances in communication, sensing, and computation have enabled large-scale networked systems in which multiple agents cooperate using only local information and neighbor-to-neighbor communication. Representative applications include economic dispatch in power networks~\cite{baranwal2020robust}, distributed estimation in sensor networks~\cite{li2007rate}, and coverage control~\cite{schwager2006distributed}. A basic problem underlying these applications is
\begin{align}
\min_{\{x_i\in\mathbb{R}^d\}_{i=1}^m}
&\ \sum_{i=1}^m f_i(x_i),
\nonumber\\[-1mm]
\text{s.t.}\quad
&x_i=x_j,\qquad i,j\in\{1,\ldots,m\},
\label{eq:DOP}
\end{align}
where agent $i$ knows only its convex local objective $f_i$ and communicates over a fixed graph.

Classical distributed gradient methods and exact decentralized refinements such as EXTRA and DIGing established convergence using local computations and communication~\cite{nedic2009distributed,shi2015extra,nedic2017achieving}. Accelerated decentralized methods subsequently incorporated momentum, gradient tracking, primal--dual updates, or repeated consensus operations~\cite{jakovetic2014fast,qu2019accelerated,xu2020OPTRA,li2024AccGT}. These developments reveal an important distinction between computational and communication complexity: acceleration with respect to gradient evaluations may require multiple communication rounds within each iteration.

Continuous-time dynamics provide a complementary perspective on acceleration. They permit Lyapunov, variational, and energy-based arguments that can expose structures that are difficult to identify directly in discrete time~\cite{budhraja2022breaking,garg2023accelerating}. Nevertheless, an accelerated differential equation does not automatically yield a rate-preserving iterative method. Straightforward Euler-type discretizations may destroy the continuous-time energy structure, and even stable discretizations need not preserve the rate of the original flow~\cite{suh2022continuous,zhang2018direct}. Understanding which guarantees survive discretization is therefore a central question.

A second issue concerns the measure of distributed suboptimality. Many analyses evaluate the centralized objective at a network-averaged or consensual surrogate~\cite{jakovetic2014fast,qu2019accelerated,li2024AccGT}. In contrast, problem~\eqref{eq:DOP} naturally leads to the aggregate objective evaluated at the actual local iterates. This quantity is sensitive to disagreement at first order and can be negative away from consensus. We therefore study its absolute deviation from the optimum together with the consensus residual.

Our main contributions are summarized below.

\noindent\textbf{1. Accelerated primal--dual flow and conserved energy.}
We introduce a second-order primal--dual flow whose time-dilated dynamics admit an Euler--Lagrange representation. A completed scaled Hamiltonian is exactly conserved and yields an $\mathcal O(t^{-2})$ rate for the aggregate objective gap and the squared consensus error. The same analysis gives an $\mathcal O(t^{-2})$ agent-wise objective rate and an $\mathcal O(t^{-1})$ disagreement rate.

\noindent\textbf{2. A barrier for single-loop discretizations.}
For a broad class of single-loop finite-memory primal--dual methods with one gradient and one communication update per iteration, we construct a horizon-dependent smooth convex instance with an $\Omega(k^{-1})$ aggregate-objective gap. Hence, this class cannot satisfy a $\mathcal O(k^{-2})$ worst-case guarantee.

\noindent\textbf{3. Rate-preserving double-loop discretization.}
Motivated by the lower bound, we combine finite-step polynomial consensus with an accelerated outer update. The resulting distributed method uses one local gradient evaluation and at most $m-1$ communication rounds per outer iteration, maintains exact consensus, and achieves an $\mathcal O(k^{-2})$ aggregate-objective rate under the same assumptions as the flow.
We compare the proposed method with DGD, EXTRA, DIGing, D-NC~\cite{jakovetic2014fast}, Acc-DNGD~\cite{qu2019accelerated}, AccGT+CA~\cite{li2024AccGT}, and OPTRA~\cite{xu2020OPTRA}. The experiments use common initialization, separate tuning and evaluation instances, and report both gradient-evaluation and communication complexity.

\section{Preliminaries}
\label{sec:Prelim}

\noindent\textbf{Notation:}
Let $N$ be a positive integer, and consider $m\geq1$ agents. For a differentiable function $g:\mathbb{R}^N\to\mathbb{R}$, its gradient is denoted by $\nabla g$. Unless specified otherwise, $\|v\|$ denotes the Euclidean norm, $\langle\cdot,\cdot\rangle$ denotes the Euclidean inner product, and $\otimes$ denotes the Kronecker product. We use $I_N$ and $\mathbf{1}_N$ for the $N\times N$ identity matrix and the $N$-dimensional all-ones vector, respectively.

Let $a_{ij}=a_{ji}\geq0$ denote the communication weight between agents $i$ and $j$ for the undirected graph. Let $\mathcal N_i\coloneqq\{j:a_{ij}>0\}$ denote the neighbor set of agent $i$. Let $\cL$ denote the corresponding graph Laplacian. Its smallest nonzero eigenvalue is denoted by $\underline{\lambda}_{\cL}$. Define $\tilde{\cL}\coloneqq\cL\otimes I_d$ and $\mathsf J\coloneqq m^{-1}\mathbf{1}_m\mathbf{1}_m^\top\otimes I_d$, 

For $X=[x_1^\top,\ldots,x_m^\top]^\top\in\mathbb{R}^{md}$, let $\bar x\coloneqq m^{-1}\sum_{i=1}^m x_i$ and $\bar X\coloneqq\mathsf JX=\mathbf{1}_m\otimes\bar x$. Define $f(x)\coloneqq\sum_{i=1}^m f_i(x)$, $f_\star\coloneqq\min_{x\in\mathbb{R}^d}f(x)$, and $F(X)\coloneqq\sum_{i=1}^m f_i(x_i)$. For any optimal consensus value $x_\star$ of~\eqref{eq:DOP}, let $X_\star\coloneqq\mathbf{1}_m\otimes x_\star$.

We use the following standard definitions.

\begin{definition}
A differentiable function $g:\mathbb{R}^N\to\mathbb{R}$ is convex if $\nabla g(x)^\top(y-x)\leq g(y)-g(x)$ for all $x,y\in\mathbb{R}^N$.
\end{definition}

\begin{definition}
A differentiable function $g:\mathbb{R}^N\to\mathbb{R}$ is $L_g$-smooth if $\|\nabla g(x)-\nabla g(y)\|\leq L_g\|x-y\|$ for all $x,y\in\mathbb{R}^N$.
\end{definition}

\begin{definition}
The Bregman divergence generated by a differentiable function $g:\mathbb{R}^N\to\mathbb{R}$ is $D_g(x,y)\coloneqq g(x)-g(y)-\langle\nabla g(y),x-y\rangle$. If $g$ is convex, then $D_g(x,y)\geq0$.
\end{definition}

Our analysis relies on the following assumptions.

\begin{assumption}
\label{assump_1}
For every agent $i\in\{1,\ldots,m\}$, the local objective $f_i\in C^1$ is convex. Moreover, $f_\star> -\infty$ and $\arg\min_{x\in\mathbb{R}^d} f(x)\neq\emptyset$.
\end{assumption}

\begin{assumption}
\label{assump_2}
The communication graph is fixed, undirected, and connected.
\end{assumption}

\begin{assumption}
\label{assump_3}
The cumulative cost function $F$ is $L_f$-smooth for some $L_f>0$.
\end{assumption}

Since $F$ is block separable, Assumption~\ref{assump_3} implies that each local objective has an $L_f$-Lipschitz continuous gradient. Separate local smoothness constants are then unnecessary.

\section{Primal--Dual Accelerated Gradient Flow and Conserved Energy}
\label{sec:AGM-ODE}

Continuous-time models such as the accelerated-gradient and optimized-gradient ODEs provide useful energy-based interpretations of acceleration for centralized convex optimization~\cite{su2014differential,suh2022continuous}. We now introduce a second-order primal--dual flow whose variational structure leads to an exact conservation law. Throughout this section, Assumptions~\ref{assump_1}--\ref{assump_3} are in force.

By Assumption~\ref{assump_2},
$\ker(\tilde{\cL})=\operatorname{span}\{\mathbf{1}_m\}\otimes\mathbb{R}^d$. Moreover, Assumption~\ref{assump_1} gives $(\mathbf{1}_m^\top\otimes I_d)\nabla F(X_\star)=\sum_{i=1}^m \nabla f_i(x_\star)=0$.
Since $\tilde{\cL}$ is symmetric, this implies
$\nabla F(X_\star)\in\operatorname{range}(\tilde{\cL})$.
Consequently, there exists a dual vector $\Lambda_\star$ satisfying the
primal--dual optimality
\begin{equation}
    \tilde{\cL}X_\star=0,
    \qquad
    \nabla F(X_\star)+\tilde{\cL}\Lambda_\star=0.
    \label{eq:pd-optimality}
\end{equation}

\subsection{Primal--dual accelerated flow}

For $t\ge t_0>0$, consider the coupled dynamics
\begin{subequations}\label{eq:HBDist}
\begin{align}
 \ddot X+\frac{3}{t}\dot X+\nabla F(X)
 +\tilde{\cL}X+\tilde{\cL}\Lambda
 +\frac{t}{2}\tilde{\cL}\dot\Lambda &=0,
 \label{eq:X-ODE}\\
 \ddot\Lambda+\frac{3}{t}\dot\Lambda
 -\tilde{\cL}X-\frac{t}{2}\tilde{\cL}\dot X &=0.
 \label{eq:Lambda-ODE}
\end{align}
\end{subequations}
Here $X=[x_1^\top,\ldots,x_m^\top]^\top$ and
$\Lambda=[\lambda_1^\top,\ldots,\lambda_m^\top]^\top$.
Note that~\eqref{eq:HBDist} is distributed. In particular, the $i^{\text{th}}$ agent implements
\begin{align*}
     &\ddot x_i+\frac{3}{t}\dot x_i+\nabla f_i(x_i) \\
     &+\sum_{j\in\mathcal N_i}a_{ij}
     \left(x_i-x_j+\lambda_i-\lambda_j
     +\frac{t}{2}(\dot\lambda_i-\dot\lambda_j)\right)=0,\\
     &\ddot\lambda_i+\frac{3}{t}\dot\lambda_i
     -\sum_{j\in\mathcal N_i}a_{ij}
     \left(x_i-x_j+\frac{t}{2}(\dot x_i-\dot x_j)\right)=0.
\end{align*}
We introduce the time-dilated coordinates
\begin{equation}
    W(t)\coloneqq t^2\bigl(X(t)-X_\star\bigr),
    \qquad
    S(t)\coloneqq t^2\bigl(\Lambda(t)-\Lambda_\star\bigr).
    \label{eq:dilated-coordinates}
\end{equation}
Using~\eqref{eq:pd-optimality}, direct differentiation of
\eqref{eq:dilated-coordinates} transforms~\eqref{eq:HBDist} into
\begin{subequations}\label{eq:WS-ODE}
\begin{align}
     &\frac{\ddot W}{t^2}\!-\!\frac{\dot W}{t^3}
     \!+\!\nabla F\!\left(\frac{W}{t^2}\!+\!X_\star\right)\!-\!\nabla F(X_\star) \!+\!\frac{\tilde{\cL}W}{t^2}
     \!+\!\frac{\tilde{\cL}\dot S}{2t} =0,
     \label{eq:W-ODE}\\
     &\frac{\ddot S}{t^2}-\frac{\dot S}{t^3}
     -\frac{1}{2t}\tilde{\cL}\dot W =0.
     \label{eq:S-ODE}
\end{align}
\end{subequations}

\subsection{Euler--Lagrange representation}

Define the kinetic and potential energy-like terms
\begin{align*}
    &T(W,S,\dot W,\dot S,t)
    \coloneqq
    \frac{1}{2t}\bigl(\|\dot W\|^2+\|\dot S\|^2\bigr) \nonumber\\
    &+\frac14\langle\dot W,\tilde{\cL}S\rangle
    -\frac14\langle W,\tilde{\cL}\dot S\rangle,\\
    &U(W,t)
    \coloneqq
    t^3D_F\!\left(\frac{W}{t^2}+X_\star,X_\star\right)
    +\frac{1}{2t}\langle\tilde{\cL}W,W\rangle,
\end{align*}
and let
\begin{equation}
    \mathscr L(W,S,\dot W,\dot S,t)\coloneqq T-U.
    \label{eq:pd-lagrangian}
\end{equation}
The following lemma recovers the dilated flow~\eqref{eq:W-ODE}--\eqref{eq:S-ODE} from  $\mathscr L$ using Euler--Lagrange equations.

\begin{lemma}
\label{lem:EL-recovery}
Under Assumptions~\ref{assump_1}--\ref{assump_3}, the Euler--Lagrange equations associated with~\eqref{eq:pd-lagrangian}
coincide exactly with~\eqref{eq:W-ODE}--\eqref{eq:S-ODE}.
\end{lemma}

\begin{proof}
Set
\[
    X=\frac{W}{t^2}+X_\star,
    \qquad
    G(X)\coloneqq\nabla F(X)-\nabla F(X_\star).
\]
Since $\tilde{\cL}$ is symmetric,
\begin{align*}
 \frac{\partial\mathscr L}{\partial\dot W}
 &=\frac{1}{t}\dot W+\frac14\tilde{\cL}S,
 &
 \frac{\partial\mathscr L}{\partial W}
 &=-\frac14\tilde{\cL}\dot S
   -tG(X)-\frac{1}{t}\tilde{\cL}W,\\
 \frac{\partial\mathscr L}{\partial\dot S}
 &=\frac{1}{t}\dot S-\frac14\tilde{\cL}W,
 &
 \frac{\partial\mathscr L}{\partial S}
 &=\frac14\tilde{\cL}\dot W.
\end{align*}
Therefore,
\begin{align*}
 \frac{d}{dt}\frac{\partial\mathscr L}{\partial\dot W}
 -\frac{\partial\mathscr L}{\partial W}
 &=
 \frac{\ddot W}{t}-\frac{\dot W}{t^2}
 +tG(X)+\frac{1}{t}\tilde{\cL}W
 +\frac12\tilde{\cL}\dot S,\\
 \frac{d}{dt}\frac{\partial\mathscr L}{\partial\dot S}
 -\frac{\partial\mathscr L}{\partial S}
 &=%
 \frac{\ddot S}{t}-\frac{\dot S}{t^2}
 -\frac12\tilde{\cL}\dot W.
\end{align*}
Equating these expressions to zero and dividing by $t$ gives precisely
\eqref{eq:W-ODE} and~\eqref{eq:S-ODE}.
\end{proof}

\subsection{Completed scaled Hamiltonian and conservation law}

The momenta generated by~\eqref{eq:pd-lagrangian} are
\[
 p_W=\frac{1}{t}\dot W+\frac14\tilde{\cL}S,
 \qquad
 p_S=\frac{1}{t}\dot S-\frac14\tilde{\cL}W.
\]
The corresponding Hamiltonian-like quantity is
\begin{align}
    &H(t)
    \coloneqq \langle\dot W,p_W\rangle
    +\langle\dot S,p_S\rangle-\mathscr  L \nonumber\\
    &=\frac{1}{2t}\bigl(\|\dot W\|^2+\|\dot S\|^2\bigr)
    +t^3D_F(X,X_\star)
    +\frac{1}{2t}\langle\tilde{\cL}W,W\rangle.
    \label{eq:Hamiltonian}
\end{align}
Thus the scaled Hamiltonian
\begin{equation*}
    V(t)\!\coloneqq\!\frac{H(t)}{t}
    \!=\!\frac{\|\dot W\|^2+\|\dot S\|^2}{2t^2}
    +t^2D_F(X,X_\star)
    +\frac{\langle\tilde{\cL}W,W\rangle}{2t^2},
\end{equation*}
is exactly the non-integral part of the energy below.

\begin{lemma}
\label{lem:energy-conservation}
Under Assumptions~\ref{assump_1}--\ref{assump_3}, along every solution of~\eqref{eq:WS-ODE}, the quantity
\begin{align}
    E(t)
    \coloneqq{}&
    \frac{\|\dot W(t)\|^2+\|\dot S(t)\|^2}{2t^2}
    +t^2D_F(X(t),X_\star) \nonumber\\
    &+\int_{t_0}^{t}2sD_F(X_\star,X(s))\,ds
    +\frac{1}{2t^2}\langle\tilde{\cL}W(t),W(t)\rangle \nonumber\\
    &+\int_{t_0}^{t}\frac{1}{s^3}
    \langle\tilde{\cL}W(s),W(s)\rangle\,ds
    \label{eq:conserved-energy}
\end{align}
is conserved, i.e., $E(t)=E(t_0)\eqqcolon E_0$ for every $t\geq t_0$.
\end{lemma}

\begin{proof}
For a time-dependent Lagrangian, every Euler--Lagrange trajectory satisfies
$\dot H(t)=-\frac{\partial\mathscr L}{\partial t}$. 
Let
$d\coloneqq X-X_\star=W/t^2$ and $G(X)=\nabla F(X)-\nabla F(X_\star)$.
Since $\partial_tX=-2W/t^3=-2d/t$ when $W$ is held fixed,
\begin{align}
    \frac{\partial\mathscr L}{\partial t}
    ={}&-\frac{\|\dot W\|^2+\|\dot S\|^2}{2t^2}
    -3t^2D_F(X,X_\star)
    +2t^2\langle G(X),d\rangle \nonumber\\
    &+\frac{1}{2t^2}\langle\tilde{\cL}W,W\rangle.
    \label{eq:partial-time-L}
\end{align}
Differentiating $V=H/t$ and using~\eqref{eq:Hamiltonian} and
\eqref{eq:partial-time-L} yields
\begin{align*}
    \dot V(t)
    &= -\frac{1}{t}\frac{\partial\mathscr L}{\partial t}
    -\frac{H(t)}{t^2} \nonumber\\
    &=2t\bigl(D_F(X,X_\star)-\langle G(X),d\rangle\bigr)
    -\frac{1}{t^3}\langle\tilde{\cL}W,W\rangle.
\end{align*}
By the definition of the Bregman divergence and the expression for $G(X)$,
\begin{align*}
 &D_F(X,X_\star)-\langle G(X),X-X_\star\rangle \\
 &=F(X)-F(X_\star)-\langle\nabla F(X),X-X_\star\rangle=-D_F(X_\star,X).
\end{align*}
Therefore,
\begin{equation}
    \dot V(t)
    =-2tD_F(X_\star,X(t))
    -\frac{1}{t^3}\langle\tilde{\cL}W(t),W(t)\rangle.
    \label{eq:perfect-derivative}
\end{equation}
Differentiating~\eqref{eq:conserved-energy} and applying the fundamental theorem of calculus shows that its two integral terms cancel the
two terms on the right-hand side of~\eqref{eq:perfect-derivative}. Hence
$\dot E(t)=0$.
\end{proof}

Since $F$ is convex and $\tilde{\cL}\succeq0$, every term in
\eqref{eq:conserved-energy} is nonnegative. In addition, multiplying
\eqref{eq:S-ODE} by $t$ gives
\[
\frac{d}{dt}\left(
\frac{\dot S}{t}-\frac12\tilde{\cL}W
\right)=0.
\]
Hence
\begin{equation}
    \frac{\dot S(t)}{t}-\frac12\tilde{\cL}W(t)
    =
    C_S,
    \,
    C_S\coloneqq
    \frac{\dot S(t_0)}{t_0}
    -\frac12\tilde{\cL}W(t_0).
    \label{eq:S-first-integral}
\end{equation}
The conservation law and~\eqref{eq:S-first-integral} directly yield the main convergence result for~\eqref{eq:HBDist}.

\begin{theorem}
\label{thm:MainThm}
Under Assumptions~\ref{assump_1}--\ref{assump_3}, let $(X,\Lambda)$ be any solution of~\eqref{eq:HBDist}. Then, for every
$t\geq t_0$,
\begin{align}
\left|
\sum_{i=1}^m f_i(x_i(t))-f_\star
\right|
&\leq
\frac{
E_0+
2\|\Lambda_\star\|
\bigl(\sqrt{2E_0}+\|C_S\|\bigr)
}{t^2},
\label{eq:aggregate-objective-rate}\\
X(t)^\top\tilde{\cL}X(t)
&\leq \frac{2E_0}{t^2}.
\label{eq:squared-consensus-rate}
\end{align}
\end{theorem}

\begin{proof}
By Lemma~\ref{lem:energy-conservation}, $E(t)=E_0$. Since every term in
\eqref{eq:conserved-energy} is nonnegative,
\begin{equation}
    D_F(X(t),X_\star)\leq \frac{E_0}{t^2}.
    \label{eq:bregman-intermediate}
\end{equation}
The kinetic term in~\eqref{eq:conserved-energy} also gives
$\|\dot S(t)\|/t\leq\sqrt{2E_0}$. Using~\eqref{eq:S-first-integral}, $\tilde{\cL}W(t)=2\frac{\dot S(t)}{t}-2C_S$, and therefore $\|\tilde{\cL}W(t)\|
    \leq
    2\bigl(\sqrt{2E_0}+\|C_S\|\bigr)$. Since $W=t^2(X-X_\star)$ and $\tilde{\cL}X_\star=0$, we have
$\tilde{\cL}W(t)=t^2\tilde{\cL}X(t)$. Consequently,
\begin{equation}
    \|\tilde{\cL}X(t)\|
    \leq
    \frac{2\bigl(\sqrt{2E_0}+\|C_S\|\bigr)}{t^2}.
    \label{eq:linear-consensus-residual}
\end{equation}

By the definition of the Bregman divergence and the optimality condition
$\nabla F(X_\star)+\tilde{\cL}\Lambda_\star=0$,
\begin{align*}
D_F(X,X_\star)
&=
F(X)-F(X_\star)
-\langle\nabla F(X_\star),X-X_\star\rangle\\
&=
F(X)-F(X_\star)
+\langle\Lambda_\star,\tilde{\cL}X\rangle.
\end{align*}
Thus, $F(X)-F(X_\star)
=
D_F(X,X_\star)
-\langle\Lambda_\star,\tilde{\cL}X\rangle$.

\noindent Since $F(X)=\sum_{i=1}^m f_i(x_i)$ and $F(X_\star)=f_\star$, the
Cauchy--Schwarz inequality, \eqref{eq:bregman-intermediate}, and
\eqref{eq:linear-consensus-residual} yield
\begin{align*}
\left|
\sum_{i=1}^m f_i(x_i(t))-f_\star
\right|
&\leq
D_F(X(t),X_\star)
+\|\Lambda_\star\|\,\|\tilde{\cL}X(t)\|\\
&\leq
\frac{
E_0+
2\|\Lambda_\star\|
\bigl(\sqrt{2E_0}+\|C_S\|\bigr)
}{t^2},
\end{align*}
which proves~\eqref{eq:aggregate-objective-rate}.
Finally, since $W=t^2(X-X_\star)$ and $\tilde{\cL}X_\star=0$,
\[
E_0
\geq
\frac{1}{2t^2}\langle\tilde{\cL}W,W\rangle
=
\frac{t^2}{2}X^\top\tilde{\cL}X.
\]
It follows that $X(t)^\top\tilde{\cL}X(t)
\leq
\frac{2E_0}{t^2}$, which proves~\eqref{eq:squared-consensus-rate}.
\end{proof}

Theorem~\ref{thm:MainThm} establishes accelerated convergence of the aggregate objective gap and the squared consensus error. The following corollary translates these into agent-wise objective convergence and an explicit rate for the disagreement norm.

\begin{corollary}
\label{cor:agentwise-and-disagreement-rates}
Under the same conditions of Theorem~\ref{thm:MainThm},
\begin{align}
f(x_i(t))-f_\star
&\leq
\frac{2E_0}{t^2}
\left(
1+\frac{mL_f}{\underline{\lambda}_{\cL}}
\right),
\label{eq:cor-agentwise-objective-rate}\\
\|X(t)-\bar X(t)\|^2
&\leq
\frac{2E_0}{\underline{\lambda}_{\cL}t^2}.
\label{eq:cor-disagreement-rate}
\end{align}
\end{corollary}

\begin{proof}
Under Assumption~\ref{assump_2},
\[
X^\top\tilde{\cL}X
=
\left\langle
\tilde{\cL}(X-\bar X),X-\bar X
\right\rangle
\geq
\underline{\lambda}_{\cL}\|X-\bar X\|^2.
\]
Combining this inequality with~\eqref{eq:squared-consensus-rate} gives
\eqref{eq:cor-disagreement-rate}.

Fix $i\in\{1,\ldots,m\}$ and define
$X_i^{\mathrm c}\coloneqq\mathbf{1}_m\otimes x_i$. Since $X_\star=\mathbf{1}_m\otimes x_\star$ and
$\sum_{j=1}^m\nabla f_j(x_\star)=0$,
\[
D_F(X_i^{\mathrm c},X_\star)=f(x_i)-f_\star.
\]
The Bregman three-point identity, $L_f$-smoothness, and Young's inequality
yield
\[
f(x_i)-f_\star
\leq
2D_F(X,X_\star)
+
L_f\|X_i^{\mathrm c}-X\|^2.
\]
Moreover, $\|X_i^{\mathrm c}-X\|^2
=
\sum_{j=1}^m\|x_i-x_j\|^2
\leq
m\|X-\bar X\|^2$. Hence, using~\eqref{eq:bregman-intermediate} and
\eqref{eq:cor-disagreement-rate}, we obtain~\eqref{eq:cor-agentwise-objective-rate}.
\end{proof}

\section{A Single-Loop Discretization Barrier}
\label{sec:discretization-barrier}
\begin{figure*}
    \centering
    \begin{subfigure}[t]{0.32\textwidth}
        \centering
        \includegraphics[width=\linewidth]
        {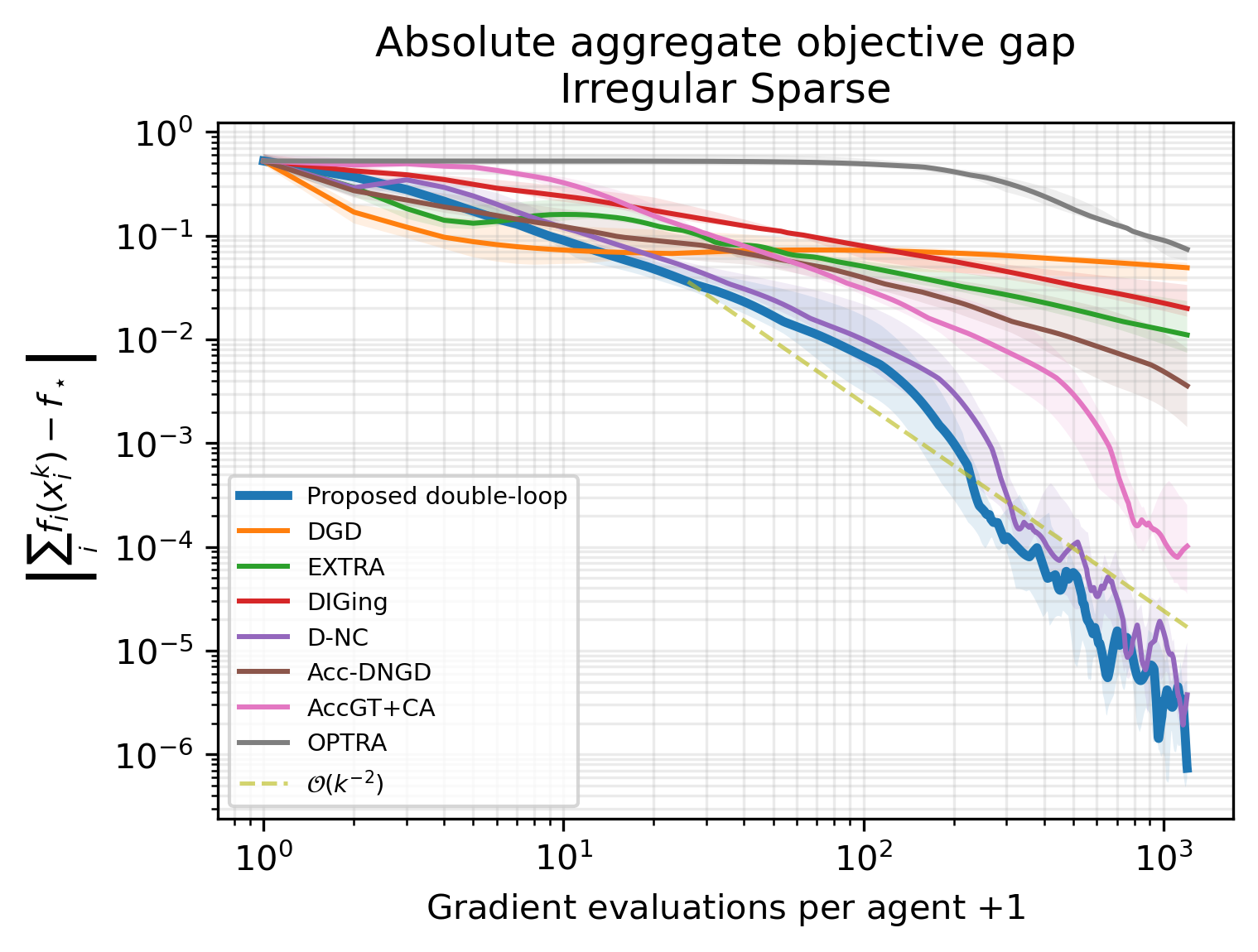}
        \caption{}
        \label{fig:objective-gradients}
    \end{subfigure}
    \hfill
    \begin{subfigure}[t]{0.32\textwidth}
        \centering
        \includegraphics[width=\linewidth]
        {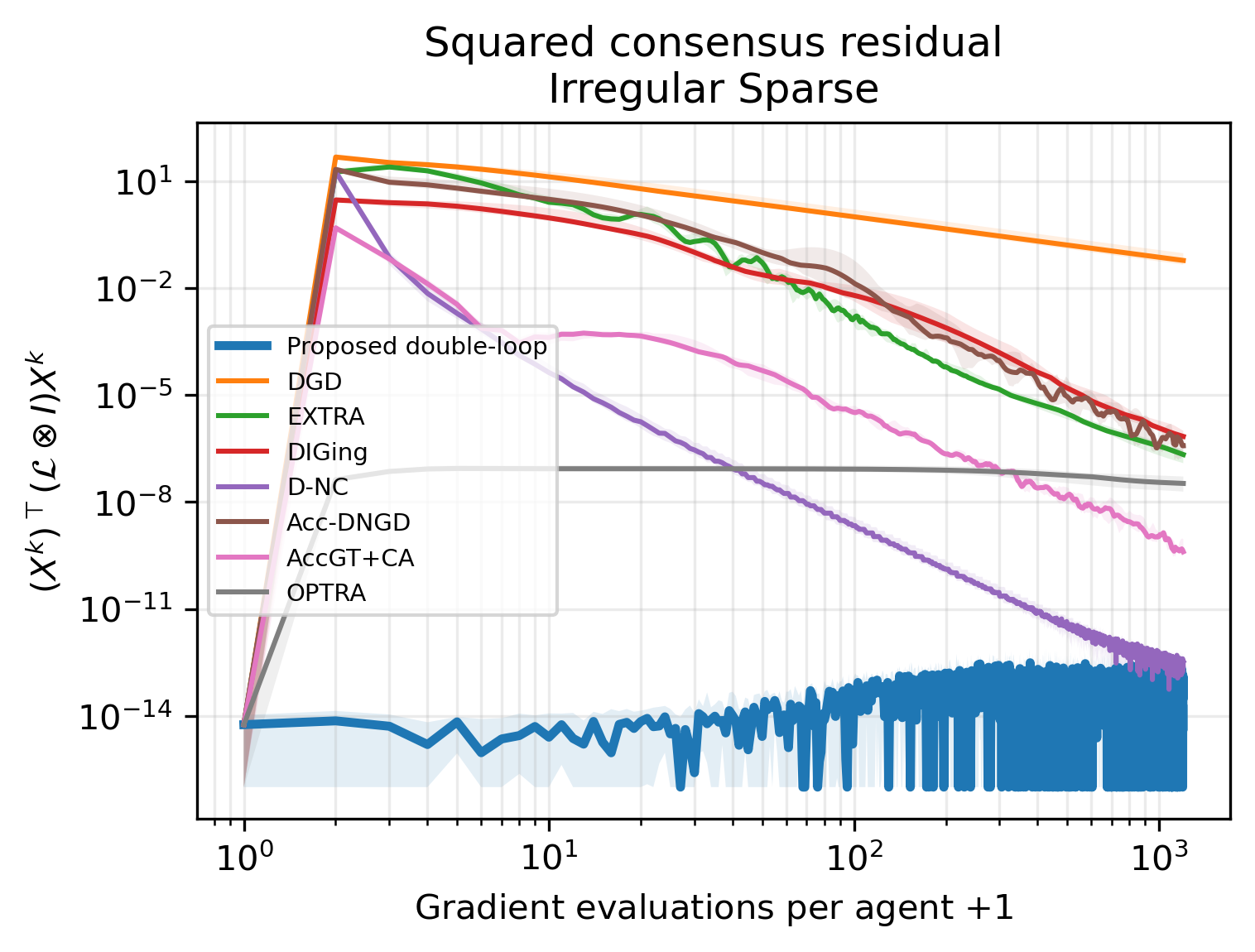}
        \caption{}
        \label{fig:consensus-gradients}
    \end{subfigure}
    \hfill
    \begin{subfigure}[t]{0.32\textwidth}
        \centering
        \includegraphics[width=\linewidth]
        {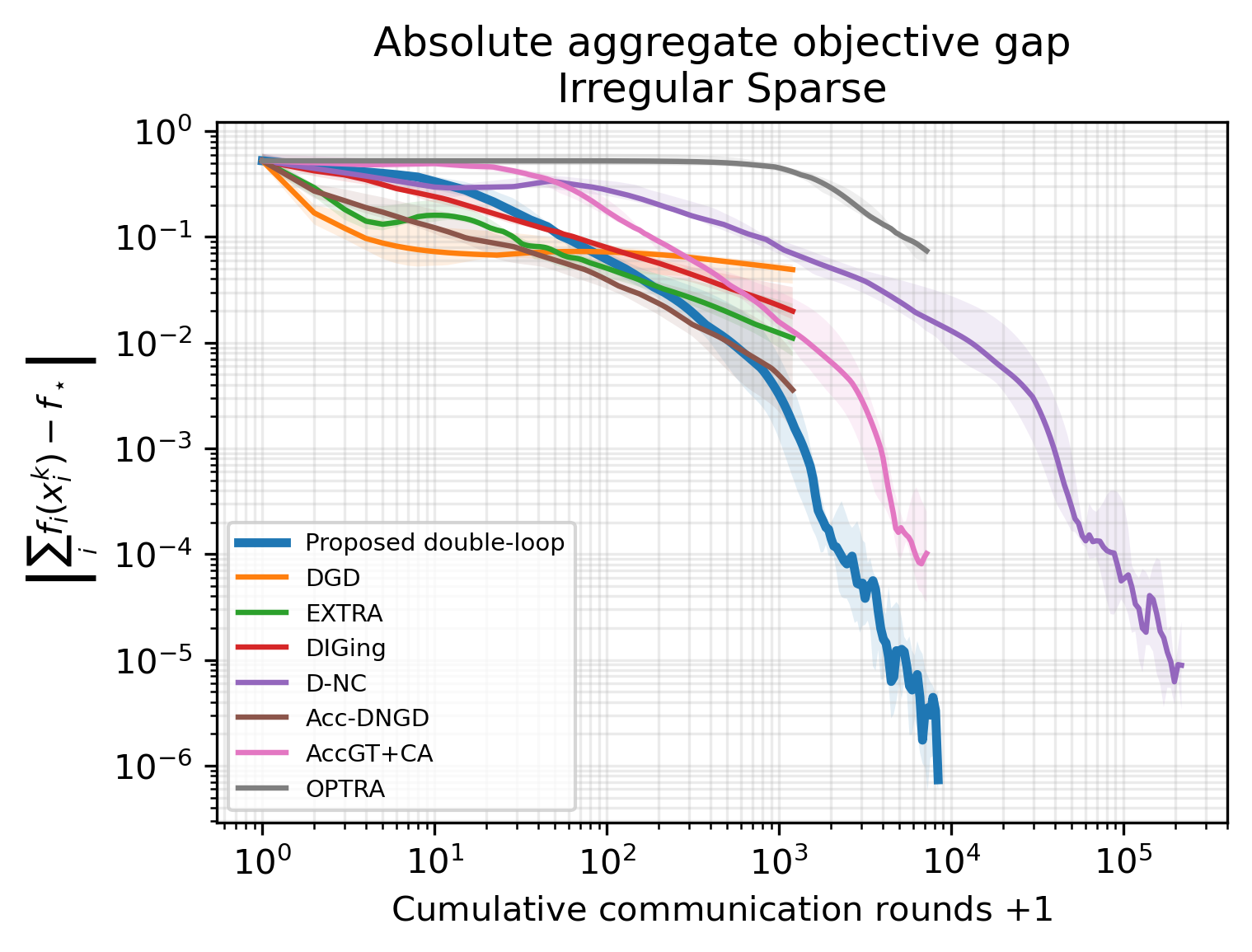}
        \caption{}
        \label{fig:objective-communications}
    \end{subfigure}
    \caption{Comparison on an irregular sparse network. The objective metric
    is $\left|\sum_{i=1}^m f_i(x_i^k)-f_\star\right|$, and the consensus
    metric is $(X^k)^\top\tilde{\cL}X^k$. Solid curves show the
    median over ten trials, and shaded regions show the interquartile range.}
    \label{fig:numerical-comparison}
    \vspace{-1.5em}
\end{figure*}

The dynamics in Section~\ref{sec:AGM-ODE} achieve an
$\mathcal{O}(t^{-2})$ rate for the aggregate objective gap. We next examine
whether this rate is inherited by a broad class of direct, single-loop
primal--dual discretizations. Fix finite memory lengths $p,q\geq0$. At iteration $k$, consider the extrapolated variables
\[
Y^k\!=\!\sum_{r=0}^{p}\delta_{r,k}X^{k-r}, 
\widehat X^k\!=\!\sum_{r=0}^{p}\gamma_{r,k}X^{k-r}, 
\widehat\Lambda^k\!=\!\sum_{r=0}^{q}\nu_{r,k}\Lambda^{k-r},
\]
where the coefficients may vary with $k$ and satisfy
$\delta_{r,k},\gamma_{r,k},\nu_{r,k}\geq0$ and
$
\sum_{r=0}^{p}\delta_{r,k}
=
\sum_{r=0}^{p}\gamma_{r,k}
=
\sum_{r=0}^{q}\nu_{r,k}
=1$.
Let $\mathsf P\in\mathbb{R}^{m\times m}$ be a doubly stochastic mixing matrix,
set $\mathsf P_d\coloneqq\mathsf P\otimes I_d$, and let
$\mathsf B\in\mathbb{R}^{md\times md}$ be any linear consensus operator satisfying $\mathsf B(\mathbf{1}_m\otimes z)=0$ for every $z\in\mathbb{R}^d$. We consider the
single-loop primal--dual recursion
\begin{subequations}\label{eq:general-discrete-class}
\begin{align}
X^{k+1}
&=
\mathsf P_dY^k
-\alpha_k\left(
\nabla F(Y^k)+\mathsf B^\top\widehat\Lambda^k
\right),
\label{eq:general-primal-update}\\
\Lambda^{k+1}
&=
\widehat\Lambda^k+\beta_k\mathsf B\widehat X^k,
\label{eq:general-dual-update}
\end{align}
\end{subequations}
where $0\leq\alpha_k\leq1/L_f$ and $\beta_k\geq0$.

This template permits time-varying stepsizes, arbitrary finite primal and
dual memory, and moving-average extrapolation, while using only one gradient
evaluation and one primal--dual communication per iteration. The
memoryless choice $p=q=0$ contains standard primal--dual gradient and
Arrow--Hurwicz-type distributed updates. For a suitable $\tau>0$, the choice $\mathsf P_d=I_{md}-\tau\tilde{\cL}$ recovers consensus-corrected variants, while nonzero $p$ and $q$ include delayed, averaged, and finite-memory versions of these methods.

The following result proves that this class cannot inherit the accelerated
aggregate-objective guarantee of the continuous-time flow over
smooth convex problems.

\begin{theorem}
\label{thm:discrete-impossibility}
Consider any algorithm of the form~\eqref{eq:general-discrete-class}. For
every horizon $N\geq1$ and every $D>0$, there exists a one-dimensional
homogeneous instance $f_1=\cdots=f_m=g_N$ satisfying
Assumptions~\ref{assump_1}--\ref{assump_3}, with
$|x_i^0-x_\star|=D$, such that
\begin{equation}
\left|
\sum_{i=1}^{m}f_i(x_i^N)-f_\star
\right|
\geq
\frac{mL_fD^2}{8(N+1)}.
\label{eq:discrete-lower-bound}
\end{equation}
\end{theorem}

\begin{proof}
Fix $N\geq1$ and define $\eta_N\coloneqq D/[2(N+1)]$. Consider the shifted
Huber function
\begin{equation}
g_N(x)\coloneqq
\begin{cases}
L_f\eta_N(D-x)-\dfrac{L_f}{2}\eta_N^2,
& x\leq D-\eta_N,\\[1ex]
\dfrac{L_f}{2}(x-D)^2,
& |x-D|\leq\eta_N,\\[1ex]
L_f\eta_N(x-D)-\dfrac{L_f}{2}\eta_N^2,
& x\geq D+\eta_N.
\end{cases}
\label{eq:shifted-Huber}
\end{equation}
The function $g_N$ is convex, continuously differentiable, and
$L_f$-smooth, with unique minimizer $x_\star=D$ and minimum value zero.
Set $f_i=g_N$ for every agent and initialize
$X^{-p}=\cdots=X^0=\mathbf{0}_m$ and
$\Lambda^{-q}=\cdots=\Lambda^0=\mathbf{0}_{md}$.

We first observe that the iterates remain consensual. Indeed, if all primal
histories are of the form $X^{k-r}=\mathbf{1}_m z_{k-r}$ and the dual histories
are zero, then $Y^k=\mathbf{1}_m y_k$ and $\widehat X^k=\mathbf{1}_m\widehat z_k$. Since $\mathsf B\mathbf{1}_m=0$, the dual update gives
$\Lambda^{k+1}=0$. Moreover,
$\mathsf P\mathbf{1}_m=\mathbf{1}_m$ and the local objectives are identical, so
the primal update gives $X^{k+1}=\mathbf{1}_m z_{k+1}$.

Define the scalar deficit $u_k\coloneqq D-z_k$. As long as the extrapolated
point lies in the left linear region of~\eqref{eq:shifted-Huber}, we have
$g_N'(y_k)=-L_f\eta_N$, and hence
\[
u_{k+1}
=
\sum_{r=0}^{p}\delta_{r,k}u_{k-r}
-\alpha_kL_f\eta_N.
\]
Let $\mu_k\coloneqq\min\{u_k,\ldots,u_{k-p}\}$. Since the extrapolation
coefficients form a convex combination,
$\mu_{k+1}\geq\mu_k-\alpha_kL_f\eta_N$. Using
$\alpha_kL_f\leq1$ and the initial value $\mu_0=D$ gives
\[
\mu_k\geq D-k\eta_N\geq D-N\eta_N>\frac{D}{2},
\qquad 0\leq k\leq N.
\]
Since $\eta_N\leq D/2$, this estimate also verifies inductively that all
points visited up to iteration $N$ remain in the left linear region. In particular, $u_N\geq D/2$, and therefore
\[
g_N(z_N)
=
L_f\eta_Nu_N-\frac{L_f}{2}\eta_N^2
\geq
\frac{L_f\eta_ND}{4}
=
\frac{L_fD^2}{8(N+1)}.
\]
All agents have the same iterate and $f_\star=0$. Thus,
\[
\left|
\sum_{i=1}^{m}f_i(x_i^N)-f_\star
\right|
=
m g_N(z_N)
\geq
\frac{mL_fD^2}{8(N+1)},
\]
which proves~\eqref{eq:discrete-lower-bound}.
\end{proof}

Theorem~\ref{thm:discrete-impossibility} is a worst-case lower bound. For each horizon $N$, it constructs a smooth convex instance on which
the aggregate objective gap after $N$ iterations is at least of order
$1/N$. Consequently, no algorithm in~\eqref{eq:general-discrete-class} admits a
worst-case $\mathcal{O}(N^{-2})$ aggregate-objective guarantee over
the smooth convex problem class. Thus, the $\mathcal{O}(t^{-2})$ aggregate-objective rate of the ODE system does not transfer to the broad
single-loop finite-memory discrete algorithm class
\eqref{eq:general-discrete-class}.

The obstruction already occurs for scalar and identical local objectives,
consensual primal initialization, zero dual initialization, and trajectories
that remain in exact consensus. It is therefore not caused by network
disagreement, objective heterogeneity, or imperfect dual tracking. Rather,
the shifted-Huber instance exposes a finite-step limitation of the underlying
first-order discretization. In particular, an accelerated aggregate-objective
rate for a discrete method requires an analysis or mechanism beyond a direct
single-loop discretization of the conserved continuous-time dynamics.

\section{Rate-Preserving Double-Loop Discretization}
\label{sec:double-loop-discretization}

In this section, we develop a double-loop discretization that preserves the
$\mathcal{O}(t^{-2})$ aggregate-objective rate of the continuous-time
dynamics in Section~\ref{sec:AGM-ODE}. Although continuous-time models often
provide useful principles for constructing accelerated algorithms, a direct
discretization need not preserve the convergence rate of the underlying
flow~\cite{suh2022continuous}. Rate-preserving discretizations based on more
elaborate numerical schemes are possible~\cite{zhang2018direct}, but may
require additional regularity or computational structure. Moreover,
Theorem~\ref{thm:discrete-impossibility} shows that the desired
aggregate-objective rate cannot hold uniformly for the broad single-loop
finite-memory class considered in Section~\ref{sec:discretization-barrier}.

Motivated by this limitation, we introduce multiple communication steps
within each outer iteration. The inner loop computes the consensus
projection of the local gradients, while the outer loop discretizes the
accelerated dynamics on the consensus manifold. The resulting method uses
one local gradient evaluation and a finite number of neighbor-communication rounds
per outer iteration.

\noindent {\bf Projected discretization:}
Using the consensus projector $\mathsf J$ defined in Section~\ref{sec:Prelim}, define $\phi(x)\coloneqq m^{-1}\sum_{i=1}^m f_i(x)$. Applying $\mathsf J$ to
\eqref{eq:X-ODE} and using $\mathsf J\tilde{\cL}=0$ gives
$\mathsf J\ddot X+(3/t)\mathsf J\dot X+\mathsf J\nabla F(X)=0$. On the
consensus manifold $X=\mathbf{1}_m\otimes x$, this equation reduces to
$\ddot x+\frac{3}{t}\dot x+\nabla\phi(x)=0$. We discretize this projected accelerated flow by combining an exact inner
consensus process with an accelerated outer update.

Let
$0=\rho_1<\rho_2<\cdots<\rho_{s_{\cL}}$ denote the distinct eigenvalues of $\cL$, assumed available to the agents for the fixed graph. For any
$Z\in\mathbb{R}^{md}$, initialize $Z^{[0]}=Z$ and do
\begin{equation}
    Z^{[\ell-1]}
    =
    \left(
    I_{md}-\frac{1}{\rho_\ell}\tilde{\cL}
    \right)Z^{[\ell-2]},
    \qquad
    \ell=2,\ldots,s_{\cL}.
    \label{eq:finite-consensus-loop}
\end{equation}
By the spectral decomposition of $\cL$,
\begin{equation}
    Z^{[s_{\cL}-1]}
    =
    \prod_{\ell=2}^{s_{\cL}}
    \left(
    I_{md}-\frac{1}{\rho_\ell}\tilde{\cL}
    \right)Z
    =
    \mathsf JZ.
    \label{eq:finite-consensus-projector}
\end{equation}
Indeed, the polynomial in~\eqref{eq:finite-consensus-projector} equals one
on $\ker(\cL)$ and vanishes on every positive-eigenvalue eigenspace.
So, exact averaging is obtained after $s_{\cL}-1\leq m-1$ communications.

\noindent {\bf Double-loop algorithm:}
The agents may start from arbitrary vectors $\xi_i^0$. Before the first
outer iteration, they apply~\eqref{eq:finite-consensus-loop} to
$\Xi^0=[(\xi_1^0)^\top,\ldots,(\xi_m^0)^\top]^\top$ and set
$X^0=Y^1=\mathsf J\Xi^0=\mathbf{1}_m\otimes x^0$, $x^0=\frac{1}{m}\sum_{i=1}^m\xi_i^0$.
Let $\theta_1=1$. At each outer iteration $k\geq1$, every agent first
evaluates one local gradient and then executes the inner consensus loop.

\begin{tcolorbox}[colback=white, colframe=black, sharp corners]
For each outer iteration $k\geq1$:
\begin{enumerate}[noitemsep, leftmargin=*]
    \item Each agent evaluates $ g_i^{k,[0]}=\nabla f_i(y_i^k)$.

    \item For $\ell=2,\ldots,s_{\cL}$, each agent performs
    \[
        g_i^{k,[\ell-1]}
        \!=\!
        g_i^{k,[\ell-2]}
        \!-
        \frac{1}{\rho_\ell}
        \sum_{j\in\mathcal N_i}a_{ij}
        (
        g_i^{k,[\ell-2]}\!-g_j^{k,[\ell-2]}
        ).
    \]

    \item Each agent applies the outer accelerated updates
    \[
        x_i^k
        =
        y_i^k-\frac{1}{L_f}g_i^{k,[s_{\cL}-1]},
        \,
        \theta_{k+1}
        =
        \frac{1+\sqrt{1+4\theta_k^2}}{2},
    \]
    \[
        y_i^{k+1}
        =
        x_i^k+
        \frac{\theta_k-1}{\theta_{k+1}}
        \left(x_i^k-x_i^{k-1}\right).
    \]
\end{enumerate}
\vspace{-1em}
\end{tcolorbox}

The inner variables are reset at every outer iteration. Since
$\nabla F(Y^k)=[\nabla f_1(y_1^k)^\top,\ldots,
\nabla f_m(y_m^k)^\top]^\top$, the inner loop is distributed and requires
only exchanges between neighboring agents. Each outer iteration therefore
uses one local gradient evaluation and $s_{\cL}-1$ communication rounds.
Each message is $d$-dimensional, so the per-neighbor communication volume is $\mathcal{O}((s_{\cL}-1)d)$ per outer iteration, while the local gradient-computation cost is the same as in standard distributed first-order methods. The additional communication rounds allow the method to avoid the single-loop lower bound of Theorem~\ref{thm:discrete-impossibility}. The following result proves that the aggregate objective gap converges at rate
$\mathcal{O}(k^{-2})$, while consensus is exact at every outer iteration.

\begin{theorem}
\label{thm:double-loop-rate}
Suppose that Assumptions~\ref{assump_1}--\ref{assump_3} hold. Let
$\{X^k\}_{k\geq0}$ be generated by the double-loop algorithm above. Then
$X^k=\mathbf{1}_m\otimes x^k$ for every $k\geq0$, and, for every $k\geq1$,
\begin{equation}
    \left|
    \sum_{i=1}^m f_i(x_i^k)-f_\star
    \right|
    \leq
    \frac{2mL_f\|x^0-x_\star\|^2}{(k+1)^2}.
    \label{eq:double-loop-objective-rate}
\end{equation}
\end{theorem}

\begin{proof}
Let $G^{k,[\ell]}\coloneqq[(g_1^{k,[\ell]})^\top,\ldots,(g_m^{k,[\ell]})^\top]^\top$. Eqn.~\eqref{eq:finite-consensus-projector} implies that the output of the inner loop at iteration $k$ is
\[
    G^{k,[s_{\cL}-1]}
    =
    \mathsf J\nabla F(Y^k)
    =
    \mathbf{1}_m\otimes
    \frac{1}{m}\sum\nolimits_{i=1}^m\nabla f_i(y_i^k).
\]
The initialization is consensual. Hence, by induction,
$X^k=\mathbf{1}_m\otimes x^k$ and
$Y^k=\mathbf{1}_m\otimes y^k$ for every $k$. The outer updates therefore
reduce exactly to
\[
    x^k
    =
    y^k-\frac{1}{L_f}\nabla\phi(y^k),
    \qquad
    y^{k+1}
    =
    x^k+
    \frac{\theta_k-1}{\theta_{k+1}}
    (x^k-x^{k-1}).
\]

By Assumption~\ref{assump_1}, $\phi$ is convex. Also,
Assumption~\ref{assump_3} gives
\begin{align*}
    &\|\nabla\phi(x)-\nabla\phi(z)\| \\
    &=
    \frac{1}{m}
    \left\|
    (\mathbf{1}_m^{\top}\otimes I_d)
    \left(
    \nabla F(\mathbf{1}_m\otimes x)
    -\nabla F(\mathbf{1}_m\otimes z)
    \right)
    \right\| \\
    &\leq
    L_f\|x-z\|,
\end{align*}
and hence $\phi$ is $L_f$-smooth. Thus, the outer recursion is the
accelerated gradient method applied to $\phi$ with step size $1/L_f$.

Let $p(y)\coloneqq y-L_f^{-1}\nabla\phi(y)$. The descent inequality for a
convex $L_f$-smooth function gives, for $z\in\mathbb{R}^d$,
\[
    \phi(z)-\phi(p(y))
    \geq
    (L_f/2) \|p(y)-y\|^2
    +
    L_f\langle y-z,p(y)-y\rangle.
\]
Define $v_k\coloneqq \phi(x^k)-\phi(x_\star)$,
$u_k\coloneqq \theta_kx^k-(\theta_k-1)x^{k-1}-x_\star$, and
$d_{k+1}\coloneqq x^{k+1}-y^{k+1}$. Applying the preceding inequality
at $y=y^{k+1}$ with $z=x_\star$ and $z=x^k$, weighting the two
inequalities by $\theta_{k+1}$ and
$\theta_{k+1}(\theta_{k+1}-1)$, respectively, and adding them, we use
$\theta_{k+1}^2-\theta_{k+1}=\theta_k^2$,
$\theta_{k+1}y^{k+1}-(\theta_{k+1}-1)x^k-x_\star=u_k$, and
$u_{k+1}=u_k+\theta_{k+1}d_{k+1}$ to obtain
\[
    (2\theta_{k+1}^2/L_f) v_{k+1}+\|u_{k+1}\|^2
    \leq
    (2\theta_{k+1}^2/L_f) v_k+\|u_k\|^2.
\]
The first gradient step gives
$(2/L_f)v_1+\|u_1\|^2\leq\|x^0-x_\star\|^2$. Therefore,
$(2\theta_k^2/L_f)v_k\leq\|x^0-x_\star\|^2$. Since
$\theta_k\geq(k+1)/2$,
$\phi(x^k)-\phi(x_\star)
    \leq
    \frac{2L_f\|x^0-x_\star\|^2}{(k+1)^2}$.
Finally, exact consensus gives
$\sum_{i=1}^m f_i(x_i^k)-f_\star
    =
    m\bigl(\phi(x^k)-\phi(x_\star)\bigr)
    \geq0$.
Multiplying the preceding estimate by $m$ proves
\eqref{eq:double-loop-objective-rate}.
\end{proof}

\section{Numerical Results}
\label{sec:numerical-results}

We evaluate the proposed method on a distributed logistic regression problem over an irregular sparse graph. Each trial uses $m=8$ agents, dimension $d=30$, and $100$ heterogeneous samples per agent, and all methods are run for $1200$ outer iterations. The comparison includes DGD, EXTRA, DIGing, D-NC~\cite{jakovetic2014fast},
Acc-DNGD~\cite{qu2019accelerated}, AccGT+CA~\cite{li2024AccGT}, and
OPTRA~\cite{xu2020OPTRA}. All methods use the same consensual initialization
and problem instances. Hyperparameters are selected on separate tuning
instances and subsequently fixed for evaluation. The curves report
median over ten independent trials, while the shaded regions indicate
inter-quartile range.

Figure~\ref{fig:numerical-comparison} reports the absolute aggregate
objective gap and the squared consensus residual against the number of local
gradient evaluations. Since several methods perform multiple communication
rounds per gradient evaluation, we also report the aggregate objective gap
against cumulative communication rounds.
The first two panels evaluate the objective and feasibility quantities
appearing in the theoretical analysis. In particular, the proposed method
maintains exact consensus at every outer iteration and exhibits the predicted
accelerated decay of the aggregate objective gap. The communication-aware
comparison in Fig.~\ref{fig:objective-communications} additionally
illustrates the cost of the inner consensus loop.

\section{Conclusion}\label{sec:Conclusion}
We introduced a second-order primal--dual flow for smooth convex distributed optimization and derived an exactly conserved energy yielding $\mathcal O(t^{-2})$ aggregate-objective and squared-consensus rates. We also established a $\Omega(k^{-1})$ lower bound for a broad single-loop finite-memory discretization class, showing that the continuous-time rate does not transfer uniformly to this class. Motivated by this barrier, we developed a double-loop method based on finite-step polynomial consensus and accelerated outer updates. The method maintains exact consensus and achieves an $\mathcal O(k^{-2})$ aggregate-objective rate, with the numerical results illustrating the accompanying communication cost.

\bibliographystyle{unsrt}
\bibliography{refs}

\begin{thebibliography}{10}

\bibitem{baranwal2020robust}
Mayank Baranwal, Kunal Garg, Dimitra Panagou, and Alfred~O Hero.
\newblock Robust distributed fixed-time economic dispatch under time-varying topology.
\newblock {\em IEEE Control Systems Letters}, 5(4):1183--1188, 2020.

\bibitem{li2007rate}
Junlin Li and Ghassan AlRegib.
\newblock Rate-constrained distributed estimation in wireless sensor networks.
\newblock {\em IEEE Transactions on Signal Processing}, 55(5):1634--1643, 2007.

\bibitem{schwager2006distributed}
Mac Schwager, James McLurkin, and Daniela Rus.
\newblock Distributed coverage control with sensory feedback for networked robots.
\newblock In {\em robotics: science and systems}, pages 49--56, 2006.

\bibitem{nedic2009distributed}
Angelia Nedic and Asuman Ozdaglar.
\newblock Distributed subgradient methods for multi-agent optimization.
\newblock {\em IEEE Transactions on Automatic Control}, 54(1):48--61, 2009.

\bibitem{shi2015extra}
Wei Shi, Qing Ling, Gang Wu, and Wotao Yin.
\newblock Extra: An exact first-order algorithm for decentralized consensus optimization.
\newblock {\em SIAM Journal on Optimization}, 25(2):944--966, 2015.

\bibitem{nedic2017achieving}
Angelia Nedic, Alex Olshevsky, and Wei Shi.
\newblock Achieving geometric convergence for distributed optimization over time-varying graphs.
\newblock {\em SIAM Journal on Optimization}, 27(4):2597--2633, 2017.

\bibitem{jakovetic2014fast}
Du{\v{s}}an Jakoveti{\'c}, Joao Xavier, and Jos{\'e}~MF Moura.
\newblock Fast distributed gradient methods.
\newblock {\em IEEE Transactions on Automatic Control}, 59(5):1131--1146, 2014.

\bibitem{qu2019accelerated}
Guannan Qu and Na~Li.
\newblock Accelerated distributed {N}esterov gradient descent.
\newblock {\em IEEE Transactions on Automatic Control}, 65(6):2566--2581, 2019.

\bibitem{xu2020OPTRA}
Jinming Xu, Ye~Tian, Ying Sun, and Gesualdo Scutari.
\newblock {Accelerated Primal-Dual Algorithms for Distributed Smooth Convex Optimization over Networks}.
\newblock In {\em International Conference on Artificial Intelligence and Statistics}, pages 2381--2391. PMLR, 2020.

\bibitem{li2024AccGT}
Huan Li and Zhouchen Lin.
\newblock {Accelerated Gradient Tracking over Time-Varying Graphs for Decentralized Optimization}.
\newblock {\em Journal of Machine Learning Research}, 25(274):1--52, 2024.

\bibitem{budhraja2022breaking}
Param Budhraja, Mayank Baranwal, Kunal Garg, and Ashish Hota.
\newblock Breaking the convergence barrier: Optimization via fixed-time convergent flows.
\newblock In {\em Proceedings of the AAAI Conference on Artificial Intelligence}, volume~36, pages 6115--6122, 2022.

\bibitem{garg2023accelerating}
Kunal Garg and Mayank Baranwal.
\newblock Accelerating distributed optimization via fixed-time convergent flows.
\newblock {\em IFAC-PapersOnLine}, 56(2):1235--1240, 2023.

\bibitem{suh2022continuous}
Jaewook~J Suh, Gyumin Roh, and Ernest~K Ryu.
\newblock Continuous-time analysis of accelerated gradient methods via conservation laws in dilated coordinate systems.
\newblock In {\em International Conference on Machine Learning}, pages 20640--20667. PMLR, 2022.

\bibitem{zhang2018direct}
Jingzhao Zhang, Aryan Mokhtari, Suvrit Sra, and Ali Jadbabaie.
\newblock Direct runge-kutta discretization achieves acceleration.
\newblock {\em Advances in neural information processing systems}, 31, 2018.

\bibitem{su2014differential}
Weijie Su, Stephen Boyd, and Emmanuel Candes.
\newblock A differential equation for modeling {N}esterov’s accelerated gradient method: {T}heory and insights.
\newblock {\em Advances in neural information processing systems}, 27, 2014.

\end{thebibliography}

\end{document}